\documentclass[12pt,leqno]{amsart}
\usepackage{amsmath,stmaryrd}
\usepackage{amssymb}
\usepackage{amsthm}
\usepackage{epsf}
\usepackage{graphicx}
\usepackage{esint} 
\usepackage{dsfont}
\usepackage[pagebackref]{hyperref} 
\hypersetup{pdfpagemode=FullScreen,  colorlinks=true} 
\usepackage{enumerate} 

\numberwithin{equation}{section}

\newcommand{\ms}{\medskip}

\newcommand{\R}{\mathbb{R}}

\renewcommand{\d}{\partial}
\newcommand{\dist}{\,\mathrm{dist}}

\newcommand{\sm}{\setminus}

\newcommand{\diam}{\mathrm{diam}}

\newcommand{\wt}{\widetilde}

\newcommand{\cW}{{\mathcal  W}}

\newcommand{\cF}{{\mathcal  F}}

\newcommand{\1}{{\mathds 1}}

\DeclareMathOperator{\diver}{div}

\newcommand{\uu}{u}
\newcommand{\Dag}{D_{\beta}^{1-\gamma}}

\theoremstyle{plain}
\newtheorem{theorem}[equation]{Theorem}
\newtheorem{lemma}[equation]{Lemma}

\newtheorem{proposition}[equation]{Proposition}

\newtheorem{definition}[equation]{Definition}

\theoremstyle{definition}

\theoremstyle{remark}

\newcommand{\dv}{\operatorname{div}}

\newcommand{\re}{\mathbb{R}}

\renewcommand{\div}{\operatorname{div}}

\newcommand{\bp}{\noindent {\it Proof}.\,\,}
\newcommand{\ep}{\hfill$\Box$ \vskip 0.08in}

\begin{document}

\title[Green estimates on complements of rectifiable sets]{Green function estimates on complements of low-dimensional uniformly rectifiable sets}

\author[David]{Guy David}
\address{Guy David. 
Universit\'e Paris-Saclay, CNRS, Laboratoire de math\'ematiques d'Orsay, 
91405, Orsay, France}
\email{guy.david@universite-paris-saclay.fr}

\author[Feneuil]{Joseph Feneuil}
\address{Joseph Feneuil. Universit\'e Paris-Saclay, CNRS, Laboratoire de math\'ematiques d'Orsay, 
91405, Orsay, France}
\email{joseph.feneuil@gmail.com}

\author[Mayboroda]{Svitlana Mayboroda}
\address{Svitlana Mayboroda. School of Mathematics, University of Minnesota, Minneapolis, MN 55455, USA}
\email{svitlana@math.umn.edu}
\thanks{
David was partially supported by the European Community H2020 grant GHAIA 777822, and the Simons Foundation grant 601941, GD. 
Mayboroda was supported in part by the NSF grant DMS 1839077 and the Simons foundation grant 563916, SM}

\begin{abstract}
It has been recently established in \cite{DMgreencomp} that on uniformly rectifiable sets the Green  function is almost affine in the weak sense, and moreover, in some scenarios such Green function estimates are equivalent to the uniform rectifiability of a set. The present paper tackles a strong analogue of these results, starting with the ``flagship" degenerate operators on sets with lower dimensional boundaries. 

We consider the elliptic operators $L_{\beta,\gamma} =- \dv D^{d+1+\gamma-n} \nabla$
associated  to a domain $\Omega \subset \R^n$ with a uniformly rectifiable boundary $\Gamma$
of dimension $d < n-1$, the now usual distance to the boundary $D = D_\beta$
given by $D_\beta(X)^{-\beta} = \int_{\Gamma} |X-y|^{-d-\beta} d\sigma(y)$
for $X \in \Omega$, where $\beta >0$ and $\gamma \in (-1,1)$. In this paper we show that the  Green  function $G$ for $L_{\beta,\gamma}$, with pole at infinity,
is well approximated by multiples of $D^{1-\gamma}$, in the sense that the function
$\big| D\nabla\big(\ln\big( \frac{G}{D^{1-\gamma}} \big)\big)\big|^2$ satisfies a Carleson measure estimate on $\Omega$. We underline that the strong and the weak results are different in nature and, of course, at the level of the proofs: the latter extensively used compactness arguments, while the present paper relies on some intricate integration by parts and the properties of the ``magical" distance function from \cite{DEM}.
 
\ms

\noindent
{\sc R\'esum\'e.} 
Dans \cite{DMgreencomp} il est d\'emontr\'e que pour les domaines \`a bord uniform\'ement rectifiable, la fonction de Green v\'erifie des estimations faibles de bonne approximation par des fonctions affines, avec une r\'eciproque vraie dans certains cas encourageants.
Ici on part de la rectifiabilit\'e uniforme et on d\'emontre les estimations fortes naturelles
d'approximation de la fonction de Green, et aussi des solutions, par des applications affines 
(ou, de mani\`ere \'equivalente, des multiples de la distance au bord adoucie). L'\'etude inclut les
analogues naturels du Laplacien dans les domaine dont la fronti\`ere est de 
grande co-dimension.

On consid\`ere les op\'erateurs elliptiques 
$L_{\beta,\gamma} = \dv D^{d+1+\gamma-n} \nabla$ associ\'es \`a un domaine
$\Omega \subset \R^n$ dont le bord $\Gamma$ est Ahlfors r\'egulier et uniform\'ement rectifiable
de dimension  $d < n-1$ et \`a la distance au bord maintenant usuelle $D = D_\beta$ d\'efinie par
$D_\beta(X)^{-\beta} = \int_{\Gamma} |X-y|^{-d-\beta} d\sigma(y)$
pour $X \in \Omega$, o\`u $\beta >0$ et $\gamma \in (-1,1)$ sont des param\`etres
et $\sigma$ une mesure Ahlfors r\'eguli\`ere sur $\Gamma$. Les auteurs ont montr\'e
pr\'ec\'edemment que la mesure elliptique associ\'ee \`a $L_{\beta,\gamma}$ est
bien d\'efinie et est mutuellement absolument continue par rapport \`a $\sigma$,
avec un poids de $A_\infty$. Ici on d\'emontre que  la fonction de Green $G$ avec p\^ole
\`a l'infini associ\'ee \`a $L_{\beta,\gamma}$ est bien approch\'ee par les multiples de 
$D$, au sens o\`u la fonction $\big| D\nabla\big(\ln\big( \frac{G}{D^{1-\gamma}} \big)\big)\big|^2$ v\'erifie une condition de Carleson sur $\Omega$. Ces nouvelles estimations sont diff\'erentes en nature.
Les estimations de \cite{DMgreencomp} reposaient sur des arguments de compacit\'e; ici on a 
besoin d'estimations plus pr\'ecises, obtenues par int\'egration par parties et en utilisant les 
propri\'et\'es alg\'ebriques de la fonction $D_\alpha$ dans le cas``magique'' de \cite{DEM}.
\end{abstract}

\maketitle

\ms\noindent{\bf Key words/Mots cl\'es.}
Uniformly rectifiable sets, degenerate elliptic operators, Estimates on Green functions, harmonic measure in higher codimension.

\ms\noindent

\tableofcontents

\section{Introduction}

Rectifiable sets are an important notion in geometric measure theory 
and the calculus of variation, in particular 
because the sets that minimize an energy often enter this category. 
In the past decades, many mathematicians worked on finding  
characterizations of rectifiability by properties apparently unrelated to geometric measure theory. 
In the early 90's, the quantifiable version - uniform rectifiability -
was introduced in \cite{DS1,DS2}
along with many characterizations in terms of geometry (such as big pieces of Lipschitz images, or using Peter Jones' $\beta$ numbers) and in terms of singular integrals.
Later, it was observed that uniformly
rectifiable sets may be the right extension of Lipschitz graphs 
for elliptic boundary value problems, that is, 
if $\Omega\subset \R^n$ is an open set with uniformly 
rectifiable boundary and  
$\Omega$ provides enough access its the boundary, then we can control the oscillations of the harmonic functions in $\Omega$. It is even more noteworthy that a criterion for rectifiability can be obtained using harmonic functions. 
Indeed, Hofmann, Martell, and Uriarte-Tuero proved in \cite{HM1,HMU} that under some conditions regarding the  
access to the boundary,  
$\partial \Omega$ is uniformly rectifiable if and only if the harmonic measure 
on $\d\Omega$
is absolutely continuous in a quantitative way - 
called 
$A_\infty$ - with respect to the surface measure 
(see also \cite{AHMNT}).  
The optimal topological conditions in this regard have been identified in \cite{AHMMT}. 
These were accompanied but a rich array of beautiful and difficult alternative characterizations, exploring Carleson estimates for the solutions, behavior of the singular integral operators, extensions to more general elliptic operators, to  mention just a few. Here  we do not aim to provide a survey of the related literature; the reader can consult, e.g., \cite{DFMKenig} for a more detailed
presentation of the literature.

A weakness of the above theory is the fact that the harmonic measure on $\Gamma$ only makes sense when $\Gamma \subset \R^n$ is of dimension $d> n-2$, because lower
dimensional sets have probability zero to be hit by a Brownian motion. 
As a consequence, rectifiability, a notion that exists for all integer dimensions, can only be characterized 
by means of the harmonic measure 
for sets $\Gamma \subset \R^n$ of dimension $d=n-1$. To overcome this obstacle, the authors of the present article developed a theory of degenerate 
elliptic operators. The idea 
was to define a `harmonic measure' on a set $\Gamma$ of dimension $d < n-1$ by replacing the Laplacian by an operator on $\Omega :=\R^n \sm \Gamma$ in the form $L = -\div w(x) \nabla$, where $w(x)$ 
goes to infinity at an appropriate rate when $x$ is approaching $\Gamma$, 
so that the corresponding `Brownian motion' is attracted by $\Gamma$ and hits 
the boundary with probability one.
The articles \cite{DFMprelim,DFMprelim2} 
set the elliptic theory for this. In particular an appropriate elliptic measure is constructed for a large
class of sets of any dimension $d<n$ (or even mixed dimensions) and operators $L$ as above, with the usual nice properties
such as the non-degeneracy and doubling properties for the harmonic measure, the Harnack inequality and
the comparison principle for solutions, and estimates for the change of poles.

Then we tested the relevancy of our new elliptic measure. 
Dahlberg proved in \cite{Da} that the classical harmonic measure is absolutely 
continuous with respect to the surface measure, 
and even $A_\infty$, whenever the boundary is Lipschitz. 
In \cite{DFMAinfty} we assumed that $\Gamma$ is
the graph of a Lipschitz function $\varphi :\R^d \to \R^{n-d}$ with small Lipschitz constant,
and proved the same $A_\infty$ property for correctly chosen operators. Yet
we had to be careful about the operator we picked because as showed in \cite{DFMKenig}, based on counterexamples from \cite{CFK,MM}, not every operator $L$ will work. 
We wanted an explicit operator, which could be constructed by a single method for any set of dimension $d$. We set our choice on
\begin{equation} \label{defLbg}
L_{\beta,\gamma}:= - \div (D_\beta)^{d+1+\gamma-n} \nabla,
\end{equation}
where $\beta >0$, $\gamma \in (-1,1)$, and $D_\beta$ is defined on $\Omega$ as
\begin{equation} \label{defDb}
D_\beta (X) := \left( \int_\Gamma |X-y|^{-d-\beta} d\sigma(y) \right)^{-1/\beta}.
\end{equation}
The quantity $D_\beta$ is equivalent to the distance to the boundary, 
i.e.,
there exists $C>0$ such that
\begin{equation} \label{Dbdistintro}
C^{-1} \dist(X,\Gamma) \leq D_\beta \leq C \dist(X,\Gamma) \qquad \text{ for } X\in \Omega,
\end{equation}
but the advantage of $D_\beta$ over $\dist(.,\Gamma)$ is to always 
be smooth in a certain quantitative way. 

When $\gamma=0$, the level of the degeneracy $D_\beta^{d+1-n}$ in the coefficients of 
the operator $L_{\beta,0}:= - \div (D_\beta)^{d+1-n} \nabla$ makes it a perfect analogue 
of the Laplacian for the sets with a $d$-dimensional boundary when $d<n-1$. 
In particular, we proved that the harmonic measure associated to $L_{\beta,\gamma}$ is $A_\infty$ 
with respect to the $d$-dimensional Hausdorff measure.
Moreover, recently, in two distinct papers \cite{DM,F} we 
extended this result to the more general case where $\Gamma \subset \R^n$ is uniformly rectifiable. 
Studying these operators and  the dimension of the support 
of the corresponding elliptic measure, we were naturally drawn to introducing a parameter 
$\gamma$ which seemingly unbalances the situation. Besides, for $d=n-1$, 
the operator $L_{\beta,\gamma}:= - \div (D_\beta)^{\gamma} \nabla$ is the celebrated 
Caffarelli-Silvester extension of the fractional Laplacian operator (cf. \cite{CS}). 
However, we were surprised to realize that the argument  in \cite{DFMAinfty} extends to all 
$\gamma\in (-1,1)$ rather simply, and the generalization
to uniformly rectifiable set is stated for any $\gamma \in (-1,1)$ in \cite{F}. 
The proof in \cite{DM} relies on geometric arguments, such as corona decompositions and the construction of sawtooth domains, and an extrapolation argument, which allow one to reduce to the case of small Lipschitz graphs.
The proof of \cite{F} is substantially simpler and more direct, and
is based on a trick unique to the case where the dimension of $\Gamma$ is at most $n-2$.

\begin{theorem}[\cite{DM,F}] \label{Main1}
Let $\Gamma \subset \R^n$ be a $d$-Ahlfors regular uniformly rectifiable set with $d < n-1$, 
and let $\sigma$ be an Ahlfors regular measure that satisfies \eqref{defAR}. 
Take $\beta >0$, $\gamma \in (-1,1)$,
define $L_{\beta,\gamma}$ as in \eqref{defLbg}, and construct the associated harmonic measure 
$\omega_{\beta,\gamma}^X$ as in Definition \ref{defhm}.
Then $\omega_{\beta,\gamma}^X$ is 
$A_\infty$-absolutely continuous 
with respect to $\sigma$.
This means that for every choice of $\epsilon \in (0,1)$, 
there exists $\delta \in (0,1)$, that depends only on $C_\sigma$, $C_0$, $\epsilon$, $n$, $d$, $\beta$, and $\gamma$,
such that for each choice of $x\in \Gamma$, $r>0$, a Borel set $E\subset B(x,r) \cap \Gamma$, and a corkscrew point $X = A_{x,r}$ as in \eqref{corkscrew}, 
\begin{equation} \label{AiTh4}
\frac{\sigma(E)}{\sigma(B(x,r) \cap \Gamma)} < \delta 
\Rightarrow  \frac{\omega^X_{\Omega,L}(E)}{\omega^X_{\Omega,L}(B(x,r) \cap \Gamma)} < \epsilon.
\end{equation}
\end{theorem}

It is known that in the present context where all our measures are doubling, the $A_\infty$ condition also implies that under the assumptions of \eqref{AiTh4},
\begin{equation} \label{AiTh}
\frac{\omega^X_{\Omega,L}(E)}{\omega^X_{\Omega,L}(B(x,r) \cap \Gamma)} < \delta \Rightarrow \frac{\sigma(E)}{\sigma(B(x,r) \cap \Gamma)} < \epsilon.
\end{equation}

Observe that in the previous theorem, contrary to the case of co-dimension 1, 
i.e. when $d=n-1$, we don't assume any topological condition. And that is perfectly natural, since the domain $\Omega:=\R^n \setminus \Gamma$ has a lot of paths and ample access to the boundary (see Lemmas~2.1 and 11.6 in \cite{DFMprelim}). The next big objective 
would be to prove
the reverse implication, meaning that 
if the harmonic measure on $\Gamma$ is $A_\infty$ with respect to the Hausdorff measure, then
$\Gamma$ 
is uniformly
rectifiable. Unfortunately (and surprisingly)
this fails brutally when $d+\beta = n+2$ (see \cite{DEM}), and although we expect this to be
the only exception, the methods used  to prove the converse in co-dimension 1 do 
not appear to be adaptable to the higher codimension case. 

The purpose of the present paper is to provide different estimates on the harmonic functions, 
which we hope 
will ultimately furnish one side of the desired criterion. Indeed, as established in \cite{DMgreencomp}, some weak bounds on the Green function are equivalent to the uniform rectifiability even in lower dimensional settings.  The general idea is that instead of trying to characterize rectifiable sets 
using the harmonic measure,  we would do so using the property 
that the Green function behaves like a distance to the boundary. 
This is not a surprising idea, as  the Green functions and harmonic measure are deeply connected, and are both prominent in the analysis of the free  boundary problems, and the proof of the properties of the harmonic measure in \cite{DFMprelim,DFMprelim2} heavily relies on a comparison between the harmonic measure and Green functions. 
Yet, the results we are about to  prove 
are  not known in the ``classical" setting of domains with  an
$n-1$ dimensional 
boundary, except for a few simple situations and certainly
not in the generality of the uniformly rectifiable sets which we attack in this paper.

We will use the Green function with a pole at infinity, which 
is constructed in \cite[Definition 6.2, Lemma 6.5]{DEM} with the following properties. 

\begin{proposition}[\cite{DEM}] \label{defGinfty}
Let $\Gamma \subset \R^n$ be a $d$-Ahlfors regular uniformly rectifiable set with $d < n-1$. Take $\beta >0$, $\gamma \in (-1,1)$, and construct $L_{\beta,\gamma}$ on $\Omega :=\R^n \setminus \Gamma$ as in \eqref{defLbg}. There exists a continuous function $G^\infty = G_{\beta,\gamma}^\infty $ on $\R^n$ such that 
\begin{enumerate}[(i)]
\item $G^\infty > 0$ on $\Omega$,
\item $G^\infty =0$ on $\Gamma$,
\item there exists a positive Borel measure $\omega_\infty = \omega^\infty_{\beta,\gamma}$ on $\Gamma$ such that 
\[\int_\Omega \nabla G^\infty \cdot \nabla \varphi \, D_\beta^{d+1-n+\gamma} dx = \int_\Gamma \varphi(y) d\omega^\infty(y) \qquad \text{ for } \varphi \in C^\infty_0(\R^n).\]
\end{enumerate}
We call $G^\infty$ the Green function with pole at infinity, and $G^\infty$ is unique up to multiplication by a scalar constant.
\end{proposition}

The Green function with pole at infinity for the Laplacian in 
$\R^n_+ = \{(x,t) \in \R^{n-1} \times (0,+\infty)\}$ is $G^\infty(x,t) = t$. 
The one for the operator $-\div |t|^{d+1-n}$ in $\R^n \setminus \R^d 
= \{(x,t) \in \R^n \times \R^{n-d}, \, t\neq 0\}$ is $G^\infty(x,t) = |t|$. 
There is a third case where $G^\infty$ can be computed, which is when 
$\Gamma$ is a $d$-Ahlfors regular set with $d < n-2$, and we choose the specific operator
$L_{n-d-2,0}$; then 
one can easily compute
that $G^\infty = D_{n-d-2}$. 
In \cite{DEM} this is called the magic case.
The rough idea of what follows
is that, when $\gamma = 0$, a good Green function at infinity should behave  
a bit like the distance to the boundary $\Gamma$, or like $D_\beta$ which is its smooth substitute. 
When $\gamma \neq 0$, similar homogeneity considerations lead us to expect that in the good
cases, $G^\infty$ will behave like $D_\beta^{1-\gamma}$ instead. Here is our main result.

\begin{theorem} \label{Main2a}
Let $\Gamma \subset \R^n$ be a $d$-Ahlfors regular uniformly rectifiable set with $d< n-1$, and let $\sigma$ be an Ahlfors regular measure 
on $\Gamma$
that satisfies \eqref{defAR}. For $\beta >0$ and $\gamma \in (-1,1)$, 
let $L_{\beta,\gamma}$ be as in \eqref{defLbg} and
define $G^\infty_{\beta,\gamma}$ as in Proposition \ref{defGinfty}. 
Then for any $\alpha >0$, there exists $C>0$ that depends only on $C_\sigma$, $C_0$, $\alpha$, $\beta$, $\gamma$, $d$, $n$ such that for any ball $B:=B(x,r)$ centered on $\Gamma$, one has
\begin{equation} \label{GiCMintro}
\int_{B} \left|\nabla \ln\left( \frac{G^\infty_{\beta,\gamma}}{D_{\alpha}^{1-\gamma}} \right)\right|^2 \, D_\alpha^{d+2-n} \, dX \leq C \sigma(B).
\end{equation}
\end{theorem}

\noindent {\em Remarks:}
\begin{itemize}
\item Remember that $G^\infty_{\beta,\gamma}$ is only defined up to a constant. 
Yet, the above theorem makes perfect
sense, because the value of left-hand side
of \eqref{GiCMintro} does 
not change when
we replace $G^\infty_{\beta,\gamma}$ by $KG^\infty_{\beta,\gamma}$, where $K$ is a positive constant.

\item Theorem \ref{Main2a} is just the application of Theorem \ref{Main2} to the Green function  
at infinity: 
we shall generalize the estimate of Theorem \ref{Main2a} to a large class of solutions that is interesting by itself. However, we decided 
to highlight the above statement, 
which was our true purpose in our search for characterizations of rectifiability.

\item Here we decided to use the same measure $\sigma$ for the definitions of $L_{\beta,\gamma}$
and $D_{\alpha}$, but a minor modification of the proof would allow us to take two different Ahlfors regular measures $\sigma$ and $\wt\sigma$ to define $L_{\beta,\gamma}$ and and $D_{\alpha}$. 

\item As we have pointed out above, these results are not known for $d=n-1$ in the generality of the uniformly rectifiable sets. In the half-space one can see somewhat similar  estimates in \cite{FKP} and the bounds on the second derivatives of the Green function for a special class of operators  were, in disguise, obtained in \cite{HMT}. 
However, it is not clear how to deduce from \cite{HMT} a co-dimension 1 analogue of \eqref{GiCMintro}, directly or using interpolation.
\end{itemize}

In the present paper we only show that good geometric properties (the uniform rectifiability of $\Gamma$)
imply precise approximation properties of $G^\infty$ by $D_\alpha$, and the proof will rely heavily
on the $A_\infty$ property of the harmonic measure and some intermediate results in \cite{F} concerning the uniform rectifiability of $\Gamma$. We do not address the issue of the converse in this paper. 
Yet there are good reasons to believe that it may be easier to prove than for the absolute continuity 
of the harmonic measure. In a parallel paper \cite{DMgreencomp},
the authors study a less precise (weaker) approximation property of the Green function, and show that
in some cases (but where $d > n-2$) it already implies the uniform rectifiability of $\Gamma$,
while in the case of the present paper it implies another strange property of potentials defined on $\Gamma$.
But we did not manage to prove that this strange property is impossible to obtain when 
$d+\alpha \neq n+2$.

The reader should be aware that even though we think of the Green function estimate \eqref{GiCMintro} as a possible alternative to the $A_\infty$ absolute continuity of the harmonic measure for the characterization of uniform rectifiability, we already know that in the general context of elliptic operator in the form $L=-\div A\nabla$, the $A_\infty$ absolute continuity \eqref{AiTh4} - or \eqref{AiTh} - is not equivalent to \eqref{GiCMintro}.
A counterexample in $\R^n_+$, that can be extended to any codimensions using the construction (4.6) in \cite{DFMKenig}, is given in the next lines. 
Denote the running point in $\R^n_+$ as 
$(x,r) \in \R^{n-1} \times (0,+\infty)$, and observe that in this case $D_\alpha(x,r) = c_\alpha r$.
We set $b(r) := 1/(2+\cos(r))$ and we construct $L_b:= -\diver b(r) \nabla$. By uniqueness of the Green function with pole at infinity (Proposition \ref{defGinfty}), we have
\[G^\infty(x,r) = \int_0^r \frac{ds}{b(s)} = 2r + \sin(r).\]
On one hand, $r \leq G^\infty(x,r) \leq 3r$ for all $(x,r) \in \R^n_+$, and that is enough, by using a comparison principle with the harmonic measure, to show that the harmonic measure associated to $L$ is equivalent to the surface measure on $\R^{n-1}$  (see for instance Theorem 1.17 in \cite{FenGinfty} for details),  hence it is
$A_\infty$ absolutely continuous with respect to the surface measure.
On the other hand, we have
\[\begin{split}
\left|\nabla \ln\Big( \frac{G_\infty}{D_\alpha} \Big)\right| = \left| \frac{\nabla G^\infty}{G^\infty} - \frac{\nabla D_\alpha}{D_\alpha}\right| = \left| \frac{2+\cos(r)}{2r + \sin(r)} - \frac{1}{r}\right| \geq \frac1{3r}  \left|\cos(r) - \frac{\sin(r)}{r}\right|.
\end{split}\]
That is, for any ball $B(x,R) \subset \R^n$ centered on the boundary $\R^{n-1}$, we have
\[\begin{split}
\int_{B(x,R)} \left|\nabla \ln\Big( \frac{G_\infty}{D_\alpha} \Big)\right|^2 D_\alpha \, dX
\geq c_n R^{n-1} \int_0^{R/2} \left|\cos(r) - \frac{\sin(r)}{r}\right|^2 \frac{dr}r \geq c'_n R^{n-1} \ln(R),
\end{split}\]
when $R$ is large (e.g. $R\geq 100$), and where $c_n,c'_n$ are constants that depends only on $n$. As a consequence, \eqref{GiCMintro} fails for the elliptic operator $L_b$. 
This highlights the difference between Green function and harmonic measure estimates for general
elliptic operators. 

Coming back to the Green function, the reader might wonder about the comparison with the results in \cite{DMgreencomp}. It is not easy to describe the results in \cite{DMgreencomp} avoiding technicalities, but roughly speaking, we only proved there the weak statement that the set where the Green function behaves like a distance to the boundary is a Carleson-prevalent set. This carries the structural information saying that there are a lot of points and a lot of scales where the desired estimate is true, but it carries no norm control. Respectively, the methods heavily  rely on compactness arguments. In the present paper the arguments and the results are completely different, aiming at a strong norm control in the form of  \eqref{GiCMintro} and using a new idea that the ``magical" distance function identified in \cite{DEM}, being an explicit solution to a certain  PDE, can be effectively used in the integration-by-parts arguments (cf. \cite{F}). A good comparison is a familiar to many experts integration by parts with the weight $t$ which disappears when the Laplacian hits $t$, although the details in our case are necessarily considerably more involved.  

In the next section, we shall give the definitions that we skipped up to now
to lighten the introduction, such as the 
definition of 
Ahlfors regular and uniformly rectifiable sets.
We also introduce the results from \cite{F} which we will 
rely upon. 
The remainder of the article will be devoted to the proof of our main result.

\ms

We shall
use the notation $A(x)\lesssim B(x)$ when $A(x)\leq CB(x)$ and $C$ is a constant whose dependence into the various parameters will be either recalled or obvious from context. 
We also write $A(x) \approx B(x)$ when $A(x)\lesssim B(x)$ and $A(x)\gtrsim B(x)$.

\section{Definitions and anterior results} 

For the rest of the article, we take $\Gamma \subset \R^n$ and $\Omega := \R^n \setminus \Gamma$. We assume that $\Gamma$ is a $d$-Ahlfors regular set with $d<n-1$, that is $\Gamma$ is closed and there exists a measure $\sigma$ supported on $\Gamma$ and a constant $C_{\sigma} \geq 1$ such that 
\begin{equation} \label{defAR}
C_{\sigma}^{-1} r^{d} \leq \sigma (B(x,r)) \leq C_\sigma r^d \qquad \text{ for } x\in \Gamma, \ r>0.
\end{equation}
It is known that if
the above property \eqref{defAR} is true for some measure $\sigma$, then it 
is also true when $\sigma$ is replaced by $\mathcal H^d|_\Gamma$ - the $d$-dimensional Hausdorff measure restricted to 
$\Gamma$.

The Ahlfors regularity of $\Gamma$ and the low dimension $d<n-1$ are sufficient conditions to obtain the aforementioned equivalence \eqref{Dbdistintro}. Indeed, Lemma 5.1 in \cite{DFMAinfty} gives us that
\begin{equation} \label{Dbdist}
C^{-1} \dist(X,\Gamma) \leq D_\beta \leq C \dist(X,\Gamma) \qquad \text{ for } X\in \Omega,
\end{equation}
where the constant $C>0$ above depends only on $\beta >0$, $C_\sigma$ and $n-d>1$.

Moreover, Lemma 11.6 in \cite{DFMprelim} entails the existence of a constant $C$ that depends only on $C_\sigma$ and $n-d>1$ such that for any $x\in \Gamma$ and $r>0$, we can find a point $A_{x,r}$ 
such that  
\begin{equation} \label{corkscrew}
C^{-1} r \leq \dist(A_{x,r},\Gamma) \leq |A_{x,r}-x| \leq Cr.
\end{equation} 
In other words, when $\Gamma$ is $d$-Ahlfors regular with $d<n-1$, its complement automatically satisfies the interior corkscrew condition. 

\medskip  

We shall also assume that $\Gamma$ is uniformly rectifiable. Equivalent definitions of uniform rectifiability 
were given 
in \cite{DS1,DS2}, and the reader may use their preferred one, but since we will
only use the uniform rectifiability of $\Gamma$ via results from \cite{F} that rely
on the summability properties of Tolsa's $\alpha$-numbers, we will use these properties 
as a definition of uniform rectifiability. We need some notation first.

We denote by $\Xi$ the set of affine $d$-dimensional planes in $\R^n$. 
Each plane $P\in \Xi$ is associated with a measure $\mu_P$, 
which is the restriction to $P$ of the 
$d$-dimensional Hausdorff measure
(i.e. $\mu_P$ is the Lebesgue measure on the plane). 
A {\bf flat measure} is a measure $\mu$ that can be written $\mu = c\mu_P$ where $c$ is a positive constant and $P\in \Xi$. The set of flat measure is called $\mathcal F$.

We need 
the following variant of
Wasserstein distances to quantify the difference between two measures, and 
then 
measure how far 
$\sigma$ is from flat measures. 

\begin{definition} \label{D6.1}
For $x\in \R^n$ and $r > 0$, denote by $Lip(x,r)$ the set of $1$-Lipschitz functions 
$f$ supported in $\overline{B(x,r)}$, that is the set of functions $f : \R^n \to \R$ 
such that $f(y)=0$ for $y\in \R^n \sm B(x,r)$ and $|f(y)-f(z)|\leq |y-z|$ for $y,z\in \R^n$. 
The normalized Wasserstein distance 
in B(x,r)
between two measures $\sigma$ and $\mu$ is
\begin{equation} \label{a6.2}
\dist_{x,r}(\mu,\sigma) = r^{-d-1} \sup_{f\in Lip(x,r)} \Big|\int f d\sigma - \int f d\mu\Big|.
\end{equation}
The distance to flat measures is then 
defined by
\begin{equation} \label{a6.3}
\alpha_\sigma(x,r) = \inf_{\mu \in \cF}\dist_{x,r}(\mu,\sigma).
\end{equation}
\end{definition}

One can easily
check that 
when \eqref{defAR} holds,
the quantity $\alpha_\sigma$ is uniformly bounded, i.e. there exists a constant 
$C$ that depends 
only on $d$, $n$, and $C_\sigma$ such that 
$\alpha_\sigma(x,r) \leq C$ for $x\in \Gamma$ and $r>0$.

Let $\Gamma$ be a $d$-Ahlfors regular set, and $\sigma$ 
a measure that satisfies \eqref{defAR}. Tolsa's characterization of uniform rectifiability,
Theorem 1.2 in \cite{Tolsa09}, is as follows\footnote{Tolsa's characterization of rectifiability in \cite{Tolsa09} is given with dyadic cubes but one can easily check that our bound \eqref{defUR} is equivalent to Tolsa's one.}:
$\Gamma$ is uniformly rectifiable if and only if 
there exists a constant $C_0>0$ such that  
\begin{equation} \label{defUR}
\int_{0}^r \int_{\Gamma \cap B(x,r)} |\alpha_\sigma(y,s)|^2 \, d\sigma(y) \, \frac{ds}{s} \leq C_0 \sigma(B(x,r)) \qquad \text{ for } x\in \Gamma \text{ and } r>0.
\end{equation}

Here we will only use the fact that \eqref{defUR} holds when $\Gamma$ is uniformly
rectifiable. That is, we will only use \eqref{defUR} and do not need to know other properties of
uniformly rectifiable sets. The property \eqref{defUR} will allow
us to obtain additional estimates on the smooth distance $D_\beta$. 
The presentation of those bounds will be easier after the following definition.

\begin{definition} \label{defCM}
Let the function $f$ be defined on $\Omega$. We say that $f$
satisfies the Carleson measure condition when
$f\in L^\infty(\Omega)$
and
$|f(X)|^2 \dist(X,\Gamma)^{d-n} dX$ is a Carleson measure, that is,
\begin{equation} \label{defCMa}
\int_{B(x,r)} |f(X)|^2 \dist(X,\Gamma)^{d-n}dX \leq C \sigma(B(x,r))
\end{equation}
for 
$x\in \Gamma$ and $r>0$, with a constant $C$ that does not 
depend on $x$ or $r$.

Thus this is actually a quadratic Carleson condition. 
For short, we shall write 
$f \in CM$,
or $f\in CM(C)$ when we want to refer to the constant in \eqref{defCMa}.
\end{definition}

\noindent Due to \eqref{Dbdist}, 
we can replace $\dist(X,\Gamma)^{d-n}$ with $D_\beta^{d-n}(X)$, and even
choose $\beta$ to fit our purposes; we shall often do this without additional explanations.
We shall rely strongly on Lemma 1.24 in \cite{F}, which says the following. 

\begin{lemma} \label{Main11a}
Let $\Gamma$ be uniformly rectifiable,
so that \eqref{defAR} and \eqref{defUR} hold.
Let $\beta >0$.  Then there exist a scalar function $b$ and a vector function $\mathcal V$, both defined on $\Omega$, such that 
\begin{equation} \label{Dbdivfree} 
\int_{\Gamma} |X-y|^{-n}(X-y) d\sigma(y) = (b  \nabla D_\beta + \mathcal V) D_\beta^{d+1-n} \qquad \text{ for } X\in \Omega
\end{equation}
and a constant $C_1$ that depends only on $C_\sigma$, $C_0$, $\beta$, $n$, and $d$, such that
\begin{equation} \label{HbV2}
C_1^{-1} \leq b \leq C_1,
\end{equation}
\begin{equation} \label{HbV4}
D_\beta \nabla b \in CM(C_1),
\end{equation}
\begin{equation} \label{HbV3}
|\mathcal V| \leq C_1,
\end{equation}
and
\begin{equation} \label{HbV5}
\mathcal V \in CM(C_1).
\end{equation}
\end{lemma}

\noindent Observe that the left-hand side of \eqref{Dbdivfree} is divergence free. So if we
use \eqref{Dbdivfree} and we write the divergence free condition in weak terms,  
we obtain
\begin{equation} \label{Dbdivfree2} \int_\Omega (b  \nabla D_\beta + \mathcal V)\cdot \nabla \varphi \, D_\beta^{d+1-n} \, dX = 0  \qquad \text{ for } \varphi \in C^\infty_0(\Omega).
\end{equation}

We also need the following, which is Lemma 1.26 in \cite{F}. 
\begin{lemma} \label{Main12}
Let $\Gamma$ be uniformly rectifiable, 
i.e., assume that  \eqref{defAR} and \eqref{defUR} hold.
Let $\alpha,\beta >0$. Then $D_\alpha \nabla[D_\beta/D_\alpha]$ satisfies the Carleson measure condition with a constant that depends only on $C_\sigma$, $C_0$, $\alpha$, $\beta$, $n$, and $d$. 
\end{lemma}

\medskip
We 
are now finished
with the 
geometric background and 
turn to the elliptic theory. Pick $\beta >0$ and $\gamma \in (-1,1)$. 
We whall use the operator $L_{\beta,\gamma}$ constructed in \eqref{defLbg};
hence 
$L_{\beta,\gamma}$ enters the scope of the theory developed in \cite{DFMprelim2},
and in particular there is an elliptic measure $\omega^X$ which we shall describe now.

We first need a Hilbert space $W_\gamma$, which is the same as in \cite{DFMprelim2},
Definition 3.1, but is more easily defined as 
\begin{equation} \label{defWg}
W_\gamma = \big\{ u \in L^1_{loc}(\R^n), \|u\|_{\gamma} := \int_\Omega |\nabla u(X)|^2 \dist(X,\Gamma)^{d+1+\gamma - n} dX < +\infty \big\};
\end{equation}
the equivalence between the two definition is proved as in \cite[Lemma 3.3 and Lemma~5.21]{DFMprelim}.

Each $f \in W_\gamma$ has a trace on $\Gamma$, which lies in a corresponding Sobolev 
space $H_\gamma$ (which is equal to $H^{1/2}(\Gamma)$ when $\gamma = 0$); then 
we denote by $W_{\gamma,0}$ the set of functions in $W_\gamma$ with zero trace;
$W_{\gamma,0}$ is also the completion
of $C^\infty_0(\Omega)$ under the norm $\|.\|_{\gamma}$. 
And for any open set $E\subset \R^n$, we write that $u \in W_\gamma(E)$ 
[respectively, $u\in W_{\gamma,0}(E)$] if $u\varphi \in W_\gamma$ 
[respectively, $u\varphi \in W_{\gamma,0}$] for any $\varphi \in C^\infty_0(E)$.

Then there is a notion of weak solution for $L_{\beta,\gamma} = 0$, which the reader may 
also find in \cite{DFMprelim2}, such that in particular
\begin{equation} \label{weaksol}
\int_\Omega (\nabla u \cdot \nabla \varphi) \, D_\beta^{d+1+\gamma-n} = 0 \qquad \text{ for } \varphi \in C^\infty_0(\Omega).
\end{equation}
when $u$ is a weak solution for $L_{\beta,\gamma} = 0$ (in $\Omega$).
Here and below, we remove the variable $X$ and the integration symbol $dX$
from the notation when they are not entirely needed; unless otherwise specified, 
all our integrals on $\Omega$ will be against the Lebesgue measure $dX$.

With all this notation, the main properties of our elliptic measures $\omega^X$
are as follows.

\begin{definition} \label{defhm}
For each $X\in \Omega$, we can define a unique probability measure 
$\omega^X:=\omega^X_{\beta,\gamma}$ on $\Gamma$ with the following properties. 
For any 
$g \in C_0(\Gamma)$ (i.e., continuous function on $\Gamma$ and compactly supported),
the function $u_g$ defined as
\[ u_g(X) = \int_\Gamma g(y) d\omega^X(y)\]
is a weak solution to $L_{\beta,\gamma}$
and, if in addition $g$ lies in the Sobolev space $H_\gamma$,
then $u_g \in W_\gamma$ and the trace of $u_g$ is equal to $g$.
\end{definition}

The space $H_\gamma \cap C_0(\Gamma)$ is dense in $C_0(\Gamma)$
(with the sup norm), so the last condition is our way of solving a Dirichlet problem.
We are now
ready to state the general version of our main theorem. 

\begin{theorem} \label{Main2}
Let $\Gamma \subset \R^n$ be a $d$-Ahlfors regular uniformly rectifiable set with $d< n-1$, and let $\sigma$ be an Ahlfors regular measure 
on $\Gamma$
that satisfies \eqref{defAR}. For $\beta >0$ and $\gamma \in (-1,1)$, define $L_{\beta,\gamma}$ as in \eqref{defLbg}. 
Then for any $\alpha >0$, there exists $C>0$ that depends only on $C_\sigma$, $C_0$, $\alpha$, $\beta$, $\gamma$, $d$, $n$,
 such that for any ball $B:=B(x,r)$ centered on $\Gamma$ and any non-negative non identically
 zero weak solution $u$ of $L_{\beta,\gamma}u = 0$ in $\Omega \cap 3B$ 
 which lies in $W_{\gamma,0}(3B)$,
 one has
\begin{equation} \label{GiCM}
\int_{B} \left|\nabla \ln\left( \frac{u}{D_{\alpha}^{1-\gamma}} \right)\right|^2 \, D_\alpha^{d+2-n} \leq C \sigma(B).
\end{equation}
\end{theorem}

The proof of the Theorem 
will
use the uniform rectifiability 
of $\Gamma$
via Theorem \ref{Main1} and Lemma \ref{Main11a}. The Lemma will be used to estimate the left-hand side of \eqref{GiCM} via simple integration techniques, except for one more complicated term in the form
\[ \int  \left|\nabla \ln\left( \frac{u}{D_{\alpha}^{1-\gamma}} \right)\right| 
|\nabla \phi_{B,\epsilon}| \, D_\alpha^{d+1-n}, \]
where $\phi_{B,\epsilon}$ is a well chosen cut off function, 
which is related to the integral of the logarithm of the Poisson kernel and that will be estimated using Theorem \ref{Main1}, that is,
 the $A_\infty$ absolute continuity of the harmonic measure.

\section{Proof of Theorem \ref{Main2}}

Let us recall

\begin{lemma} \label{Main22}
Let $\Gamma \subset \R^n$ be a $d$-Ahlfors regular uniformly rectifiable set with $d < n-2$, and let $\sigma$ be an Ahlfors regular measure 
on $\Gamma$
that satisfies \eqref{defAR}. Define $L_{\beta,\gamma}$ as in \eqref{defLbg}. 
Then there exist $C>0$ and $\theta\in (0,1]$, that depend
only on $C_\sigma$, $C_0$, $\beta$, $\gamma$, $n$ and $d$, 
such that for each choice of $x\in \Gamma$, $r>0$, any Borel set $E \subset B(x,r) \cap \Gamma$, and any corkscrew point $X = A_{x,r}$ as in \eqref{corkscrew}, one has
\begin{equation} \label{AiTh2}
\frac{\omega_{\beta,\gamma}^X(E)}{\omega_{\beta,\gamma}^X(B(x,r))} \leq C \left( \frac{\sigma(E)}{\sigma(B(x,r))} \right)^\theta
 \end{equation}
 and
 \begin{equation} \label{AiTh3}
\frac{\sigma(E)}{\sigma(B(x,r))}\leq C \left(  \frac{\omega_{\beta,\gamma}^X(E)}{\omega_{\beta,\gamma}^X(B(x,r))} \right)^\theta.
 \end{equation}
\end{lemma}

Here and below, we assume implicitly that the  constant $C$ in \eqref{corkscrew} 
(the definition of corkscrew points) is chosen to depend on $C_\sigma$, $n$ and $d$ only. 
If we allow a  larger $C$ in \eqref{corkscrew}, the constants in \eqref{AiTh2} and \eqref{AiTh3} depend on 
$C$ as well.

\medskip

\bp 
The conditions \eqref{AiTh2} and \eqref{AiTh3} are another characterizations of the fact that $\omega^X \in A_\infty(\sigma)$, which is true 
by Theorem \ref{Main1}.
The fact that \eqref{AiTh2}--\eqref{AiTh3} is implied by \eqref{AiTh} can be found in \cite[Theorem 1.4.13]{KenigB} and its proof in \cite[Lemma 5]{CF}.
\ep

For our next result, we want to establish that the logarithm of the Poisson kernel, 
that is, 
$\ln(\frac{d\omega^X}{d\sigma})$, is integrable. We want a quantitative version, 
and moreover, we shall state this 
in a form that is more directly applicable when we
need it (in the proof of Proposition \ref{Main21}).

\begin{lemma} \label{Main23}
Let $\Gamma$, $\sigma$, and $L_{\beta,\gamma}$ as in Lemma \ref{Main22}.

Take $B:=B(x,r)$, 
a ball centered on $\Gamma$, 
and $X = A_{x,r}$, a corkscrew point as in \eqref{corkscrew}. 
If $\{Q_i\}_{i\in \mathcal I}$ is a finitely overlapping collection of Borel subsets of $B \cap \Gamma$, then
\begin{equation} \label{Main23a}\sum_{i\in \mathcal I} \left|\ln\left(\frac{\omega_{\beta,\gamma}^X(Q_i)}{\sigma(Q_i)} \frac{\sigma(B)}{\omega_{\beta,\gamma}^X(B)}\right) \right|\sigma(Q_i) \leq C \sigma(B),\end{equation}
where $C$ depends only on $C_\sigma$, $C_0$, $\beta$, $\gamma$, $n$, $d$, and the maximal number of overlaps in the collection $\{Q_i\}_{i\in \mathcal I}$.
\end{lemma}

\bp
We introduce for $k\in \mathbb Z$, 
\[\mathcal I_k := \left\{i\in \mathcal I, \, 2^k \leq \frac{\omega_{\beta,\gamma}^X(Q_i)}{\sigma(Q_i)} \frac{\sigma(B)}{\omega_{\beta,\gamma}^X(B)} \leq 2^{k+1} \right\}.\] 
Then we define
\[E_k := \bigcup_{i\in \mathcal I_k} Q_i.\]
Since
$\{Q_i\}$ is finitely overlapping, we have $\sum_{i\in \mathcal I_k} \sigma(Q_i) \leq C\sigma(E_k)$, and we can thus write
\begin{equation} \label{lnPka}
\sum_{i\in \mathcal I} \ln\left(\frac{\omega_{\beta,\gamma}^X(Q_i)}{\sigma(Q_i)} \frac{\sigma(B)}{\omega_{\beta,\gamma}^X(B)}\right) \sigma(Q_i) \lesssim \sum_{k\in \mathbb Z} (|k|+1) \sigma(E_k).
\end{equation}
Yet, 
since $\{Q_i\}_{i\in \mathcal I_k}$ 
is a finitely overlapping covering of $E_k$,
\[\omega_{\beta,\gamma}^X(E_k) \approx \sum_{i\in \mathcal I_k} \omega_{\beta,\gamma}^X(Q_i) \approx 2^k \frac{\omega_{\beta,\gamma}^X(B)}{\sigma(B)} \sum_{i\in \mathcal I_k} \sigma(Q_i) \approx  2^k \frac{\omega_{\beta,\gamma}^X(B)}{\sigma(B)} \sigma(E_k),\]
which means that 
\begin{equation} \label{insln2k}
\frac{\omega_{\beta,\gamma}^X(E_k)}{\sigma(E_k) }  \frac{\sigma(B)}{\omega_{\beta,\gamma}^X(B)} \approx 2^k.
\end{equation}
The use of 
\eqref{insln2k} in \eqref{AiTh2} leads to 
\[\sigma(E_k) \lesssim 2^{-k/(1-\theta)} \sigma(B)\]
while the use of \eqref{insln2k} in \eqref{AiTh3} gives 
\[\sigma(E_k) \lesssim 2^{k\theta/(1-\theta)} \sigma(B).\]
We use the first of the last two estimate to bound $\sigma(E_k)$ when $k \geq 0$ and the second one to bound $\sigma(E_k)$ when $k$ is negative. Combined with \eqref{lnPka}, we deduce
\[\begin{split}
\sum_{i\in \mathcal I} \ln\left(\frac{\omega_{\beta,\gamma}^X(Q_i)}{\sigma(Q_i)} \frac{\sigma(B)}{\omega_{\beta,\gamma}^X(B)}\right) \sigma(Q_i) & \lesssim \sum_{k \geq 0} (k+1) 2^{-\frac k{1-\theta}} \sigma(B) + \sum_{k < 0} (1-k) 2^{-\frac {k\theta}{1-\theta}} \sigma(B) \\
& \lesssim \sigma(B).
\end{split}\]
The lemma follows.
\ep

\begin{proposition} \label{Main21}
Let $\Gamma \subset \R^n$ be a $d$-Ahlfors regular uniformly rectifiable set 
with $d< n-2$, and 
let $\sigma$ be an Ahlfors regular measure that satisfies \eqref{defAR}. 
Take $\beta >0$ and $\gamma \in (-1,1)$, define $L_{\beta,\gamma}$ as in \eqref{defLbg}. 
Then for
any ball $B$ centered on $\Gamma$ and any non-negative non identically zero 
weak solution $u$ to $L_{\beta,\gamma} u = 0$ in $3B \cap \Omega$, 
with $u\in W_{\gamma,0(3B)}$, 
\begin{equation} \label{GiCMz}
\int_{B} \left|\nabla \ln\left( \frac{\uu}{D_{\beta}^{1-\gamma}} \right)\right|^2 \, D_\beta^{d+2-n} \leq C \sigma(B).
\end{equation}
with a constant $C>0$ that depends only on $C_\sigma$, $C_0$, $\beta$, $\gamma$, 
$d$,
and $n$.
\end{proposition}

Observe that Proposition \ref{Main21} is 
the  special case of 
Theorem \ref{Main2} where we take $\alpha =  \beta$. 
We shall now check that conversely, Theorem \ref{Main2} follows from 
Proposition \ref{Main21} and Lemma~\ref{Main12}.

\medskip
\noindent {\em Proof of Theorem \ref{Main2} from Proposition \ref{Main21} with the help of Lemma \ref{Main12}}. Let $\alpha>0$. Notice that
\[\begin{split}
\left|\nabla \ln\left( \frac{u}{D_{\alpha}^{1-\gamma}} \right)\right| & = \left|\nabla \left[ \ln\left( \frac{u}{D_{\beta}^{1-\gamma}} \right)  + (1-\gamma) \ln\left( \frac{D_{\beta}}{D_{\alpha}} \right) \right] \right| \\
& \leq \left|\nabla \left[ \ln\left( \frac{u}{D_{\beta}^{1-\gamma}} \right) \right] \right| + (1-\gamma) \frac{D_{\alpha}}{D_\beta} \left| \nabla \left[ \dfrac{D_\beta}{D_{\alpha}} \right]\right|.
\end{split}\]
So Proposition \ref{Main21}, Lemma \ref{Main12}, and \eqref{Dbdist} easily implies \eqref{GiCMz}.
\ep

\noindent{\em Proof of Proposition \ref{Main21}.} The proof of the lemma will be divided into 4 steps. The core step is step 4, where we use the properties of the solutions. Step 3 treats in advance the most complicated term that we met in Step 4. In this step, we will compare a solution $u$ to a Green function, then to the harmonic measure, and finally we use Lemma \ref{Main23}. In Step 2, we construct the finitely overlapping covering $\{Q_i\}_{i\in \mathcal I}$ that will be needed in order to invoke Lemma \ref{Main23}. At last, Step 1 introduces the cut-off function $\phi_{B,\epsilon}$ used to bound the left-hand side of \eqref{GiCM} and shows that $D_\beta \nabla \phi_{B,\epsilon}$ satisfies the Carleson measure condition.

\medskip

\noindent {\bf Step 1: Introduction of the cut-off function $\phi_{B,\epsilon}$.} Let $B=B(x,r)$ be a ball in $\R^n$ centered on the boundary $\Gamma$ and $\epsilon>0$ small. The proof of Theorem \ref{Main2} is a local one and thus as usual will involve cut off functions. Take $\psi \in C^\infty_0(\R)$ be such that $\psi \equiv 1$ on $[-1,1]$, $\psi$ is compactly supported in $(-2,2)$, $0\leq \psi \leq 1$, and $|\psi'|\leq 2$. We define the function $\phi_{B,\epsilon}$ on $\Omega$ by
\begin{equation} \label{defphi}
\phi_{B,\epsilon}(X) := \psi\left(\frac{\dist(X,B)}{10\dist(X,\Gamma)}\right) \psi\left(\frac{2\dist(X,B)}{r}\right) \psi\left(\frac{\epsilon}{\dist(X,\Gamma)}\right).
\end{equation}
The support of $\phi_{B,\epsilon}$ is thus contained in
\begin{equation} \label{defE0}
E_0 = \big\{X \in 2B, \, \dist(X,B) \leq 20\dist(X,\Gamma) \text{ and } \dist(X,\Gamma) 
\geq \epsilon/2 \big\}.
\end{equation}
The gradient  of $\phi_{B,\epsilon}$ comes from $3$ regions; the first one (associated
to the first cut-off) is 
\begin{equation} \label{defE1}
\begin{aligned}
\big\{X \in E_0 \, ,\, &10\dist(X,\Gamma) \leq \dist(X,B) \leq 20 \dist(X,\Gamma) \big\}
\\&
\subset E_1 : =
 \big\{X\in 2B, \, 10\dist(X,\Gamma) \leq \dist(X,B) \leq 20 \dist(X,\Gamma) \big\},
 \end{aligned}
\end{equation}
the second one is
\begin{equation} \label{defE2}
\big\{X \in E_0 \, , \, r \leq 2\dist(X,B) \leq 2r \big\}
 \subset E_2 : =
 \big\{X\in 2B, \, \, r/40 \leq \dist(X,\Gamma) \leq 2r \big\}
\end{equation}
because $\dist(X,\Gamma) \leq 2r$ when $X\in 2B$ and similarly
$\dist(X,\Gamma) \geq \frac{1}{20} \dist(X,B) \geq \frac{r}{20}$, and
the third region is contained in  
\begin{equation} \label{defE3}
E_3:= \big\{X\in 2B, \, \epsilon/2 \leq \dist(X,\Gamma) \leq \epsilon \big\}.
\end{equation}
Then it is easy to check that 
the gradient of $\phi_{B,\epsilon}$ satisfies 
\begin{equation} \label{gradphi}
|\nabla \phi_{B,\epsilon}| \leq \frac{100}{\dist(X,\Gamma)} \left[ \1_{E_1} + \1_{E_2} + \1_{E_3}\right].
\end{equation}
We claim that $\1_{E_1}$, $\1_{E_2}$, and $\1_{E_3}$ all satisfy the Carleson measure condition.
That is, we have
\begin{equation} \label{claimEi}
\int_{B(y,s)} \1_{E_j}^2  \dist(X,\Gamma)^{d-n} \leq Cs^d
\end{equation} 
for $1 \leq j \leq 3$, $y\in  \Gamma$, and $0 < s < +\infty$.
Of course the square can be removed, and it will follows from \eqref{gradphi}, 
\eqref{claimEi} (as soon as we prove it),
and \eqref{Dbdist} that
for any $\beta>0$, $y\in \Gamma$ and  $s\in (0,+\infty)$, we have
\begin{equation} \label{claimEi2}
\int_{B(y,s)} |D_\beta \nabla \phi_{B,\epsilon}| D_\beta^{d-n} + \int_{B(y,s)} |D_\beta \nabla \phi_{B,\epsilon}|^2 D_\beta^{d-n}\leq Cs^d,
\end{equation}
where the constant $C>0$ depends only on $\beta$, $n$, $d$, and $C_\sigma$.

We will prove the claim \eqref{claimEi} in the end of the next Step. 

\medskip

\noindent {\bf Step 2: Construction of the collection $\{Q_i\}_{i\in \mathcal I}$.} 
We would like to have a collection of boundary cubes $Q_i\subset \Gamma$ and the corresponding Whitney cubes $R_i\subset \Omega$ satisfying rather usual nice properties: bounded overlap, control of the size and  and the distance to the boundary for $R_i$'s, reasonable placement of the corkscrew points. There are plenty of papers that present various versions of this construction, but here, in fact, we need something simpler than usual as we do not need to control the cones or the related Harnack chains. For these reasons, let us simply carry out the construction by hands. 

We need a family of Whitney cubes $\mathcal W$, as constructed in \cite{Stein93}. 
We record
the basic
properties of $\mathcal W$ that we shall need. The collection $\mathcal W$ is the family
of maximal dyadic cubes $R\subset \Omega$ such that $20R \subset \Omega$, 
different cubes $R_i$ and $R_j$, $i\neq j$, in $\mathcal W$ have disjoint interiors (by maximality), 
their union covers  $\Omega = \R^n \sm \Gamma$, and for $Q \in \cW$
\begin{equation}
20Q \subset \Omega \text{ but } 60R\cap \Gamma \neq \emptyset.
\end{equation}
Moreover, given $R \in \mathcal W$, the number of cubes $\wt R \in \mathcal W$ such that
\begin{equation} \label{fglm4}
\dist(R,\Gamma) \approx \dist(\wt R,\Gamma) \ \text{ and } \,
\dist(R,\wt R) \lesssim \dist(R,\Gamma)
\end{equation}
is uniformly bounded by a constant that depends only on the dimensions and the constants involved in \eqref{fglm4}. Also, the cubes $\wt R \in \cW$ such that $3R \cap 3\wt R$,
all satisfy \eqref{fglm4}.

We will use the cubes of $\cW$ to cover the region $E_1\cup E_2 \cup E_3$.
Consider the subset $\cW_0 \subset \cW$ of cubes $R \in \cW$ that meet
$E_1\cup E_2 \cup E_3$, and  label these cubes with a set $I$, so that 
$\cW_0 = \{ R_i \, ; \, i\in I \}$. 
We also want to associate a boundary ball
$Q_i = \Gamma \cap B(x_i,r_i)$ to each $R_i$, $i\in I$. This is a classical thing to do,
but we shall do it by hand to get a slightly  better control.
We shall choose the  $Q_i$ so that
\begin{equation} \label{defWxra}
Q_i \subset 3B \cap \Gamma
\end{equation} 
(we can do this and this will be helpful because we shall consider solutions in $3B$),
\begin{equation} \label{defWxrb}
\text{$\{Q_i\}_{i\in \mathcal I}$ has bounded overlap,}
\end{equation} 
where the bound for the overlap depends only on $C_\sigma$, $n$, and $d$,
and 
\begin{equation} \label{defWxr}
r_i \approx \dist(R_i,\Gamma) \approx \dist(R_i,x_i),
\end{equation} 
also with constants depend only on $C_\sigma$, $n$, and $d$. 
As a consequence, we will be able to use any point $X_i \in R_i$ as a corkscrew point 
for the pair $(x_i,Cr_i)$ (see \eqref{corkscrew}).

So let us construct the $Q_i$. 
First write $I = I_1 \cup I_2 \cup I_3$, a disjoint union where 
$I_2 = \big\{ i \in I \, ; \, R_i \cap E_2 \neq \emptyset \big\}$,
then $I_3 = \big\{ i \in I \sm I_2 \, ; \, R_i \cap E_3 \neq \emptyset \big\}$,
and finally $I_1 = \big\{ i \in I \sm (I_2\cup I_3) \, ; \, R_i \cap E_1 \neq \emptyset \big\}$.

We start with $i  \in I_2$. This is the simplest case because $E_2 \subset 2B = B(x,2r)$,
and $\dist(X,\Gamma) \geq r/40$ on $E_2$. By definition of $\cW$, $I_2$ has at most
$C$ elements (see near \eqref{fglm4}), and for each $i \in I_2$ we take
$x_i = x$ and $r_i = 3r$ (and hence $Q_i = \Gamma \cap 3B$); the constraints
\eqref{defWxra}, \eqref{defWxrb}, \eqref{defWxr} for $I_2$ are easily checked.

Next consider $i\in I_3$; thus $R_i$ meets $E_3$, where 
$\varepsilon/2 \leq \dist(X,\Gamma) \leq \varepsilon$. Pick
$x_i \in \Gamma$ such that $\dist(x_i, E_3 \cap R_i) = \dist(\Gamma, E_3 \cap R_i)$,
and take $r_i = \varepsilon$. Then $Q_i  \subset 3B$ if $\varepsilon$ is small
enough, because $E_3 \subset 2B$ and $\diam(R_i) \approx\dist(R_i, \Gamma)
\approx \dist(E_3 \cap R_i, \Gamma) \approx \varepsilon$, \eqref{defWxr}
holds for the same reasons, and the $Q_i$, $i\in I_3$ have bounded overlap
by the property near \eqref{fglm4}.

Finally, for $i  \in I_1$, take $X_i \in  R_i \cap E_1$ and $x_i \in \Gamma$ such that
$ |X_i-x_i| = \dist(R_i \cap E_1,\Gamma)$, and set $Q_i = \Gamma \cap B(x_i,r_i)$,
with $r_i = |X_i-x_i|$. Since $X_i \in E_1$, the definition \eqref{defE1} yields
\begin{equation} \label{defWxr3}
\dist(X_i,B) \leq 20 \dist(X,\Gamma) = 20 |X_i-x_i| \leq 2\dist(X_i,B) \leq 2r,
\end{equation} 
and in particular $x_i\in \frac52B$. Also, $X_i \notin E_2$ (because $i \notin I_2$),
so $r_i = \dist(R_i \cap E_1,\Gamma) \leq \dist(X_i,\Gamma) \leq r/40$, and hence $Q_i \subset 3B$,
as  needed for \eqref{defWxra}. The Whitney property \eqref{defWxr} holds essentially
by definition (and because $R_i$ is a Whitney cube), so we just need to bound the overlap
of the $Q_i$.

Assume that $Q_j \cap Q_i \neq \emptyset$ for two indices $i,j\in \mathcal I_1$,
and let $X_i, X_j, r_i$, and $r_j$ be as above. Let us first check that $r_i \approx r_j$.
By \eqref{defWxr3}
\[
\dist(x_i,B) \leq |X_i-x_i| + \dist(X_i,B) \leq 21 |X_i - x_i| = 21 r_i
\]
and 
\[
\dist(x_i,B) \geq \dist(X_i,B) - |X_i-x_i| \leq  9|X_i - x_i| = 9 r_i,
\]
which can be summarized as
\begin{equation} \label{defWxr4}
9r_i \leq \dist(x_i,B) \leq 21 r_i.
\end{equation} 
Similarly,
\begin{equation} \label{defWxr5}
9r_j \leq \dist(x_j,B) \leq 21 r_j.
\end{equation} 
We can assume without loss of generality that
$\dist(x_j,B) \leq \dist(x_i,B)$, which, together with \eqref{defWxr4}-\eqref{defWxr5}, 
leads to $9r_j \leq 21r_i$. Moreover, since $Q_i \cap Q_j \neq \emptyset$, 
\[r_i + r_j \geq |x_j - x_i| \geq \dist(x_i,B) - \dist(x_j,B) \geq 9r_i - 21r_j,\]
hence 
$22r_j \geq 8r_i$. Recall that $9r_j \leq 21r_i$, so the two radii are
equivalent. The bounded overlap property follows, because
$ |x_j - x_i| \leq r_i + r_j$ when $Q_i \cap Q_j \neq \emptyset$.

This completes our construction of Whitney cubes $R_i$ and associated surface 
balls $Q_i$, with the properties \eqref{defWxra}--\eqref{defWxr}
(notice that the overlap constant for the $Q_i$, $i \in I$, is less than the sum
of the overlap constants for the $I_j$).

\ms
With this at hand, let us prove \eqref{claimEi}. This is now quite easy as for any $y\in  \Gamma$ and $0 < s < +\infty$ we have
\begin{multline} \label{claimEi-bis}
\int_{B(y,s)} \1_{E_1\cup E_2\cup E_3}^2  \dist(X,\Gamma)^{d-n} \\ \leq \sum_{R_i: \, R_i\cap B(y,s)\neq \emptyset} \int_{R_i}\dist(X,\Gamma)^{d-n} \lesssim \sum_{R_i: \, R_i\cap B(y,s)\neq \emptyset}  r_i^d \lesssim s^d,
\end{multline} 
using the bounded overlap property \eqref{defWxrb}.

\medskip

\noindent {\bf Step 3:  Estimates for the integral $S$.
} 
Take $u \in W_{\gamma,0}(3B)$, 
a weak solution to $L_{\beta,\gamma}u = 0$ in $3B$. 
Let $X_0 \in 2B$ be a corkscrew point for the boundary ball $2B \cap \Gamma$, and also choose another
corkscrew point $X_1 \in \Omega \cap 8B \setminus 4B$, so that  
$\dist(X_0,\Gamma) \geq cr$ and
$\dist(X_1,\Gamma) \geq cr$ for a small constant $c$ that depends only on 
$C_\sigma$, $d$, and $n$. 
In this Step 3, we prove that
\begin{equation} \label{claimS}
S:= \int_{E_1\cup E_2 \cup E_3} \left| \ln \left( \frac{u}{D_\beta^{1-\gamma}} \frac{D_\beta^{1-\gamma}(X_0)}{u(X_0)}\right)\right| D_\beta^{d-n} \leq C r^d,
\end{equation}
where $C$ depends only on $C_\sigma$, $C_0$, $\beta$, $\gamma$, $n$ and $d$.

We shall use the comparison principle to compare $u$ to 
the Green function, and then estimate 
the Green function in terms of  
harmonic measure. We define the Green function $g$ on $\Omega \times \Omega$ 
as in \cite[Section 11]{DFMprelim2}. The precise definition is not relevant for the present proof, and the properties that we care about are the fact that $X \to g(X,X_1)$ lies in $W_{\gamma,0}(3B)$ and is a solution to $L_{\beta,\gamma} u = 0$ in $3B$, and that $g(X,Y) = g(Y,X)$ (which is true because the operator $L_{\beta,\gamma}$ is selfadjoint).
Theorem 15.64 in \cite{DFMprelim2} (the 
comparison theorem) yields that 
\begin{equation} \label{claimS1}
\frac{u(X)}{u(X_0)} \approx \frac{g(X,X_1)}{g(X_0,X_1)} \qquad \text{ for } X\in 2B,
\end{equation}
with constants that depend only on $C_\sigma$, $C_0$, $n$, $d$, $\beta$ and $\gamma$. Actually, Theorem 15.64 in \cite{DFMprelim2} requires the solutions (in our case $u$ and $g(.,X_1)$) to be solutions in a larger 
ball $2KB \cap \Omega$, and not in only $3B\cap \Omega$. 
This condition is only needed 
because \cite{DFMprelim2} also allows sets $\Gamma$ of codimension $1$ or less,
and then we need
to ensure that we can connect every component of $2B \cap \Omega$ by 
Harnack chains that stays in $2KB$.
Here $\Gamma$ is of dimension $d < n-1$, so 
$2B\cap \Omega$ is very well connected in 
the first place (see Lemma 2.1 in \cite{DFMprelim}) ,
and 
assuming that $u$ and $g(.,X_1)$ are  
solutions in $3B$ is enough. Using the fact that $g(X,Y)$ is symmetric, we deduce that
\begin{equation} \label{claimS2}
\frac{u(X)}{u(X_0)} \approx \frac{g(X_1,X)}{g(X_1,X_0)} \qquad \text{ for } X\in 2B.
\end{equation}

Recall that $X_0$ is a corkscrew point associated to $2B$
and that if $X \in R_i$, then 
$\dist(X,\Gamma) \geq C^{-1} r_i$ and $|X-x_i| \leq C r_i$
by \eqref{defWxr}, so $X$ can be used as a corkscrew point for the boundary ball 
$Q_i$.
As a consequence, Lemma~15.28 in \cite{DFMprelim2}, 
(where here $m(B \cap \Omega)$ is the mass of $B$ for the measure 
$\dist(X,\Gamma)^{d+1+\gamma-n} dX$ on $\Omega$, with $\gamma \in (-1,1)$, so 
$m(B \cap \Omega) \approx r^{d+1+\gamma}$; see the discussion in Section 3.2 of 
\cite{DFMprelim2}), and the doubling property of harmonic measure (Lemma 15.43
of \cite{DFMprelim2}) give that
\[
g(X_1,X_0) \approx \dist(X_0,\Gamma)^{1-d-\gamma} \omega^{X_1}(3B)
\]
and for each $i\in I$ 
\[
g(X_1,X) \approx \dist(X,\Gamma)^{1-d-\gamma} \omega^{X_1}(Q_i) \qquad \text{ for } X \in R_i,
\]
where the constants depend only on $C_\sigma$, $C_0$, $n$, $d$, $\beta$ and $\gamma$. Using the equivalence \eqref{Dbdist}
and the Ahlfors regularity of $\sigma$,  we deduce that
\begin{equation} \label{claimS3}
g(X_1,X_0) \approx D_\beta(X_0)^{1-\gamma} \, \frac{\omega^{X_1}(3B)}{\sigma(3B)}
\end{equation}
and for each $i\in I$, 
\begin{equation} \label{claimS4}
g(X_1,X) \approx D_\beta(X)^{1-\gamma} \,\frac{\omega^{X_1}(Q_i)}{\sigma(Q_i)} \qquad \text{ for } X \in R_i.
\end{equation}
We gather \eqref{claimS2}, \eqref{claimS3}, and \eqref{claimS4} to obtain that, for every $i\in I$ 
 and every $X\in R_i$,
\begin{equation} \label{claimS5}
\frac{u(X)}{D_\beta^{1-\gamma}(X)} \frac{D_\beta^{1-\gamma}(X_0)}{u(X_0)} \approx \frac{\omega^{X_1}(Q_i)}{\sigma(Q_i)} \frac{\sigma(3B)}{\omega^{X_1}(3B)},
\end{equation}
where the constants depend only on $C_\sigma$, $C_0$, $n$, $d$, $\beta$ and $\gamma$. 
This
immediately implies that for $i\in I$ 
and $X\in R_i$,
\[\left| \ln \left( \frac{u(X)}{D_\beta^{1-\gamma}(X)} 
\frac{D_\beta^{1-\gamma}(X_0)}{u(X_0)}\right)\right| \leq C + \left| \ln \left( \frac{\omega^{X_1}(Q_i)}{\sigma(Q_i)} \frac{\sigma(3B)}{\omega^{X_1}(3B)} \right)\right|\]
and then
\begin{equation} \label{claimS6}
\begin{split}
S & 
\leq  
\sum_{i\in I}
\int_{R_i} \left| \ln \left( \frac{u}{D_\beta^{1-\gamma}} \frac{D_\beta^{1-\gamma}(X_0)}{u(X_0)}\right)\right| D_\beta^{d-n} \\
& \leq \sum_{i\in I} 
\left[C+ \left| \ln \left( \frac{\omega^{X_1}(Q_i)}{\sigma(Q_i)} \frac{\sigma(3B)}{\omega^{X_1}(3B)} \right)\right| \right]  \int_{R_i} D_\beta^{d-n}.
\end{split}
\end{equation}
By\eqref{Dbdist} and the fact that 
$R_i \in \cW$ is a Whitney cube with the property \eqref{defWxr}
(by construction), we have $D_\beta \approx r_i$ 
on $R_i$
and $|R_i| \lesssim r_i^n$. Hence 
\[\int_{R_i} D_\beta^{d-n} \lesssim r_i^d \approx \sigma(Q_i)\]
by \eqref{defAR}. Using this observation in \eqref{claimS6}, we infer that 
\[\begin{split}
S & \leq \sum_{i\in I} 
\left[C+ \left| \ln \left( \frac{\omega^{X_1}(Q_i)}{\sigma(Q_i)} \frac{\sigma(3B)}{\omega^{X_1}(3B)} \right)\right| \right]  \sigma(Q_i) \\
& =  C  \sum_{i\in I} 
\sigma(Q_i) +  \sum_{i\in I} 
\left| \ln \left( \frac{\omega^{X_1}(Q_i)}{\sigma(Q_i)} \frac{\sigma(3B)}{\omega^{X_1}(3B)} \right)\right|  \sigma(Q_i).
\end{split}\]
The first term of the right-hand side  
is bounded by $C\sigma(3B)$ since
the $Q_i$, $i\in I$, are contained in $\Gamma \cap 3B$ (by \eqref{defWxra})
and have bounded overlap (by \eqref{defWxrb}).
 The second term 
 is also less than 
$C\sigma(3B)$, by 
Lemma \ref{Main23}. We conclude that
\[S \lesssim \sigma(3B) \lesssim r^d\]
by \eqref{defAR}. The claim \ref{claimS} follows.
 
\medskip

\noindent {\bf Step 4: Core of the proof.} 
Let us turn to the main and last step of the proof. 
Set
\begin{equation} \label{dTT}
T := \int_{\Omega} \bigg| \nabla \ln 
\Big( \frac{\uu}{D_\beta^{1-\gamma}} \Big)
\bigg|^2 \phi_{B,\epsilon}^2 D_{\beta}^{d+2-n}.
\end{equation}
We aim to prove that
\begin{equation} \label{claimZ}
T \leq C r^d + C r^{d/2} T^{1/2},
\end{equation}
where $C$ depends only on $C_\sigma$, $C_0$, $\beta$, $\gamma$, $n$ and $d$. The solution $\uu$ and the smooth distance $D_\beta$ are uniformly bounded from above and from below by a positive constant on the support of the cut-off function $\phi_{B,\epsilon}$. Thanks to this fact, the 
quantities of both sides of \eqref{claimZ} are finite, and 
\eqref{claimZ} self-improves into 
\begin{equation} 
\int_{\Omega} \Big|  
\nabla \ln \Big(  
\frac{\uu}{D_\beta^{1-\gamma}} \Big) \Big|^2
\phi_{B,\epsilon}^2 D_{\beta}^{d+2-n} 
= T \leq C r^d,
\end{equation}
where $C>0$ depends on the same parameters as in \eqref{claimZ}. Once we are there, the left-hand side is uniformly bounded in $\epsilon$, so taking $\epsilon \to 0$ leads to the desired result \eqref{GiCMz}.

Keep in mind that 
\begin{equation} \label{nablalnGD}
\nabla \ln \bigg(  
\frac{\uu}{D_\beta^{1-\gamma}} \bigg)  
= \frac{\nabla \uu}{\uu} - \frac{\nabla \Dag}{\Dag} = \dfrac{\Dag \nabla \uu - \uu \nabla \Dag}{\Dag \uu}. 
\end{equation}

We shall 
use Lemma \ref{Main11a} to obtain the existence of a scalar function $b$ and a vector function $\mathcal V$ such that 
\begin{equation} \label{HbV1}
H_{n-d-1} := \int_\Gamma |X-y|^{-n}(X-y) d\sigma(y) = 
(b \nabla D_\beta + \mathcal V) D_{\beta}^{d+1-n} 
\end{equation}
as in \eqref{Dbdivfree}, with the bounds \eqref{HbV2}-\eqref{HbV5}.

Before we start to bound $T$, let us comment  
on $H_{n-d-1}$ and $\nabla D_\beta$. First, the vector function $H_{n-d-1}$ is smooth 
and divergence free in $\Omega$,
and the formulation of this fact in the weak sense is that for any compactly supported 
$\varphi \in W^{1,1}(\Omega)$, 
\begin{equation} 
\label{3.40a}
\int_\Omega H_{n-d-1} \cdot \nabla \varphi = 0.
\end{equation}
Set $H_\alpha(X) := \int_\Gamma |X-y|^{-d-1-\alpha}(X-y) d\sigma(y)$ for $\alpha > 0$.
Since $|H_\alpha(X)| \leq \int_\Gamma |X-y|^{-d-\alpha} d\sigma(y) = D_\alpha(X)^{-\alpha}$,
\eqref{Dbdist} implies that
\begin{equation}
\label{3.40b}
|H_\alpha| \lesssim D_\alpha^{-\alpha} \lesssim D_\beta^{-\alpha}.
\end{equation}
Moreover, observe that $\nabla (D_\beta^{-\beta})
= -(d+\beta) H_{\beta+1}$, 
directly by \eqref{defDb}; since also $\nabla (D_\beta^{-\beta})
= - \beta D_\beta^{-\beta -1} \nabla D_\beta$, we deduce from \eqref{3.40b} that
\begin{equation}
\label{3.40c}
|\nabla D_\beta| \lesssim 1.
\end{equation}

We turn to the proof of \eqref{claimZ}. 
Write $\phi$ instead of $\phi_{B,\epsilon}$ to lighten the notation. By \eqref{HbV2} and \eqref{nablalnGD}, 
\[ \begin{split}
T \lesssim \int_{\Omega} \bigg| 
\nabla \ln \bigg(  
\frac{\uu}{D_\beta^{1-\gamma}} \bigg) 
\bigg|^2  
\phi^2 b \, D_{\beta}^{d+2-n} & = \int_\Omega \frac{b\nabla \uu}{\uu} \cdot 
\bigg( 
\dfrac{\Dag \nabla \uu - \uu \nabla \Dag}{\Dag \uu} \bigg) 
\ \phi^2 D_\beta^{d+2-n}  \\
& \hspace{2cm} - \int_\Omega \frac{b\nabla \Dag}{\Dag}\cdot  \nabla \bigg[ 
\ln \bigg(  
\frac{\uu}{D_\beta^{1-\gamma}} \bigg)\bigg]  
\phi^2 D_\beta^{d+2-n} \\
& := T_1 - T_2.
\end{split} \]

Let us start with $T_2$.  We simplify the factors $D_\beta$ and we invoke the relation \eqref{HbV1} to get
\[ \begin{split}
T_{2} & = (\gamma-1) \int_\Omega \frac{b\nabla D_\beta}{D_\beta^{n-d-1}}\cdot  \nabla \bigg[ 
\ln \bigg( 
\frac{\uu}{D_\beta^{1-\gamma}} \bigg) 
\bigg]
\phi^2 \\
& = 
(\gamma - 1) 
\int_\Omega H_{n-d-1} \cdot  \nabla \bigg[ 
\ln \bigg( 
\frac{\uu}{D_\beta^{1-\gamma}} \bigg)\bigg] 
\phi^2 
- (\gamma - 1)  
\int_\Omega \mathcal V \cdot  \nabla \bigg[
\ln \bigg( 
\frac{\uu}{D_\beta^{1-\gamma}} \bigg)\bigg] 
\phi^2  D_\beta^{d+1-n}\\
& := T_{21} + T_{22}.
\end{split} \]
We use the Cauchy-Schwarz inequality to bound $T_{22}$, as follows:
\[ \begin{split}
|T_{22}| & \lesssim \Big( 
\int_{\Omega} |\mathcal V|^2 \phi^2  D_\beta^{d-n} 
\Big)^\frac12 
\bigg( 
\int_{\Omega} \bigg| 
\nabla \bigg[ 
\ln \bigg( 
\frac{\uu}{D_\beta^{1-\gamma}} 
\bigg)\bigg]\bigg|^2 
\phi^2  D_\beta^{d+2-n} 
\bigg)^\frac12 \\
& \lesssim  
T^{1/2} \left( \int_{2B} |\mathcal V|^2  D_\beta^{d-n}\right)^\frac12 
\lesssim r^{d/2} T^{1/2}
\end{split} \]
by \eqref{HbV5}. This fits with \eqref{claimZ}.
As for $T_{21}$, notice that its value will not be changed if we replace $u$ by $Ku$, where $K$ is a constant. Hence 
\[T_{21} =
(\gamma-1) 
\int_\Omega H_{n-d-1} \cdot  \nabla \bigg[ 
\ln \bigg( 
\frac{K\uu}{D_\beta^{1-\gamma}} \bigg)\bigg] 
\phi^2,\]
where $K$ is a constant to be chosen later.  
We force $\phi^2$ into the 
gradient, 
then use the fact that $H_{n-d-1}$ is divergence free (see \eqref{3.40a}), and 
obtain
\[ \begin{aligned}
T_{21} & = 
(\gamma-1) 
\int_\Omega H_{n-d-1} \cdot \nabla \Big[ 
\phi^2 \ln \Big(  
\frac{K\uu}{D_\beta^{1-\gamma}} \Big) \Big]  
- 2
(\gamma-1) \int_\Omega H_{n-d-1} \cdot  \nabla \phi \ \phi \ln 
\Big( 
\frac{K\uu} 
{D_\beta^{1-\gamma}} \Big)  
\\
& 
 =  0 - 2 (\gamma-1)
\int_\Omega H_{n-d-1} \cdot  \nabla \phi \ \phi \ln 
\Big(
\frac{K \uu}
{D_\beta^{1-\gamma}} \Big). 
\end{aligned} 
\]

Recall that $|H_{n-d-1}| \lesssim D_\beta^{d+1-n}$ by \eqref{3.40b}, 
and that $|\nabla \phi| \lesssim \1_{E_1 \cup E_2 \cup E_3} /D_\beta$ 
by \eqref{gradphi} and \eqref{Dbdist}, so 
\[
|T_{21}| \lesssim \int_{E_1 \cup E_2 \cup E_3} \bigg| 
\ln \bigg( 
\frac{K\uu}{D_\beta^{1-\gamma}} \bigg)\bigg|   
D_\beta^{d-n}.
\]
We choose $K = D_\beta^{1-\gamma}(X_0)/u(X_0)$, so that the right-hand side above is what we called $S$ 
in \eqref{claimS}. Hence $T_{21} \lesssim r^d$
by \eqref{claimS},  
as needed for \eqref{claimZ}.

We switch to the estimation of $T_1$. We want to use the fact that $\uu$ is a solution, 
and for this we write
\begin{equation} \begin{split}
T_1 & =  
\int_\Omega \nabla \uu \cdot \frac{D_\beta^{1-\gamma}}{\uu} \bigg( 
\dfrac{\Dag \nabla \uu - \uu \nabla \Dag}{\Dag \uu} \bigg) 
\ b\phi^2 D_\beta^{d+1+\gamma-n} \\
 & 
 = - 
 \int_\Omega \nabla \uu \cdot \nabla \bigg[ 
 \dfrac{\Dag}{\uu} \bigg] 
 \ b\phi^2 D_\beta^{d+1+\gamma-n} \\
 & 
 = - 
 \int_\Omega \nabla \uu \cdot \nabla \bigg[ 
 b\phi^2 \dfrac{\Dag}{\uu} \bigg]  
 \  D_\beta^{d+1+\gamma-n}
 +
 \int_\Omega \nabla \uu \cdot \nabla b \bigg(
 \phi^2 \dfrac{\Dag}{\uu} \bigg) 
 \  D_\beta^{d+1+\gamma-n} \\
 & 
 \qquad 
 +2 
 \int_\Omega \nabla \uu \cdot \nabla \phi \bigg( 
 \phi b\dfrac{\Dag}{\uu} \bigg) 
 \  D_\beta^{d+1+\gamma-n} \\
 & := T_{11} + T_{12} + T_{13}.
\end{split} \end{equation}
Notice that $T_{11}= 0$ because 
$\uu$ is a weak solution to $L_{\beta,\gamma} u = 0$ on $\Omega \cap 2B$
and $\phi^2 \dfrac{\Dag}{\uu}$ lies in $W^{1,2}_{loc}(\Omega)$ and compactly supported in $\Omega$.
The terms $T_{12}$ and $T_{13}$ are morally similar. In both case, we don't like the terms with 
$\uu$ because we don't know so
much about it, so we use \eqref{nablalnGD} to replace $\uu$  
by the nice function $\Dag$, and the difference will be controlled with the help of $T^{1/2}$, the square root of our initial integral. We start with
\begin{equation} \begin{split}
T_{13} & 
=  2 
\int_\Omega \frac{\nabla \uu}{\uu} \cdot \nabla \phi \, \phi\, 
b \,
D_\beta^{d+2-n} \\
& 
=  2
\int_\Omega \left( \frac{\nabla \uu}{\uu} - \frac{\nabla \Dag}{\Dag}\right) 
\cdot (\nabla \phi) \,  \phi \, b\, D_\beta^{d+2-n} 
- 2 (\gamma - 1) 
\int_\Omega \frac{b\nabla D_\beta}{D_\beta^{n-d-1}} \cdot (\nabla \phi) \phi \\
& := T_{131} + T_{132}.
\end{split} \end{equation}
We use \eqref{HbV2} and 
 the fact that $D_\beta \nabla \phi$ satisfies the Carleson measure property to get
\[\begin{split}
|T_{131}| & \lesssim \Big(  
\int_\Omega |\nabla \phi|^2 D_\alpha^{d+2-n} 
\Big)^{1/2} 
\bigg(
\int_\Omega \bigg|
\frac{\nabla \uu}{\uu} - \frac{\nabla \Dag}{\Dag} 
\bigg|^2 \phi^2 D_\alpha^{d+2-n} 
\bigg)^{1/2}  \\
& \lesssim r^{d/2} \bigg( 
\int_{\Omega} \bigg| 
\nabla \ln \bigg( 
\frac{u}{D_\alpha^{1-\gamma}} \bigg) \bigg|^2  
\phi^2 D_{\alpha}^{d+2-n}
\bigg)^{1/2}
\end{split}\]
by \eqref{claimEi2}, and then \eqref{nablalnGD}. 
The bound of $T_{131}$ that we just obtained appears in the right-hand side 
of \eqref{claimZ}, as desired. 
As for $T_{132}$, we invoke \eqref{3.40c}, \eqref{HbV2}, and then \eqref{claimEi2} to write
\[\begin{split}
|T_{132}| & \lesssim \int_{2B} |D_\beta \nabla \phi| D_\beta^{n-d} \lesssim r^d
\end{split}.\]

Similarly to $T_{13}$, we treat $T_{12}$ as follows: 
\begin{equation} \begin{split}
T_{12} & 
= 
\int_\Omega \frac{\nabla \uu}{\uu} \cdot \nabla b \, \phi^2 D_\beta^{d+2-n} \\
& 
= 
\int_\Omega \bigg(
\frac{\nabla \uu}{\uu} - \frac{\nabla \Dag}{\Dag} \bigg) 
\cdot \nabla b \, \phi^2 D_\beta^{d+2-n} 
- (\gamma - 1) 
\int_\Omega \nabla D_\beta \cdot \nabla b \, \phi^2 D_\beta^{d+1-n} \\
& := T_{121} + T_{122}.
\end{split} \end{equation}
Thanks to the definition \eqref{dTT} and then \eqref{HbV4},
\begin{equation} \begin{split}
|T_{121}| & \lesssim \bigg( 
\int_\Omega \bigg|
\frac{\nabla \uu}{\uu} - \frac{\nabla \Dag}{\Dag} 
\bigg|^2 \phi^2 D_\beta^{d+2-n}
\bigg)^{1/2} \Big(
\int_\Omega 
|\nabla b|^2  
\phi^2 D_\beta^{d+2-n}
\Big)^{1/2} \\
& \lesssim T^{1/2} \left( \int_{2B} |D_\beta \nabla b|^2 D_\beta^{d-n}\right)^{\frac12} \lesssim r^{d/2} T^{1/2}.
\end{split} \end{equation}
Now, we want to use \eqref{HbV1} again, so we force the function $b$ to appear 
and we place all the remaining terms in the second gradient. Then
$T_{122}$ becomes 
\[\begin{split}
T_{122} &
 = -
(\gamma - 1) \int_\Omega b\nabla D_\beta \cdot \frac{\nabla b}b \, \phi^2 D_\beta^{d+1-n} \\
& 
= -
(\gamma - 1) \int_\Omega b\nabla D_\beta \cdot \nabla [\phi^2 \ln(b)] \, D_\beta^{d+1-n}
- 2
(1-\gamma) \int_\Omega b\nabla D_\beta \cdot \nabla \phi \, \phi\ln(b) \, D_\beta^{d+1-n}
\\
&:= T_{1221} + T_{1222}.
\end{split}\]
We start with $T_{1222}$.
Since $b \approx 1$, we have that $|b \ln(b)| \lesssim 1$. Besides, $|\nabla D_\beta| \lesssim 1$ by \eqref{3.40c}. Therefore,
\[|T_{1222}| \lesssim \int_{2B} |D_\beta \nabla \phi| D_\beta^{d-n},\]
and then 
$|T_{1222}| \lesssim r^d$ by \eqref{claimEi2}.
We use \eqref{HbV1} to  write $T_{1221}$ as
\[\begin{split}
T_{1221} & = (1-\gamma) \int_\Omega H_{n-d-1}\cdot \nabla [\phi^2 \ln(b)] 
+ (\gamma - 1) \int_\Omega \mathcal V \cdot 
\nabla [\phi^2 \ln(b)] \, D_\beta^{d+1-n}
\\
& = 0 + (\gamma - 1) \int_\Omega \mathcal V \cdot 
\Big[2\phi\nabla \phi \ln(b) + \phi^2\frac{\nabla b}{b} \Big] D_{\beta}^{d+1-n}
\end{split}\]
by \eqref{3.40a}. 
Hence by \eqref{HbV2} and then the Cauchy-Schwarz inequality,
\[\begin{split}
|T_{1221}| & \lesssim \int_{2B} |\mathcal V| (|\nabla \phi| + |\nabla b|)  D_{\beta}^{d+1-n} \\
& \lesssim \left(\int_{2B} |\mathcal V|^2 D_\beta^{n-d} \right)^{\frac12}  \left( \int_{2B} (|D_\beta \nabla \phi|^2 + |D_\beta \nabla b|^2)  D_{\beta}^{d-n} \right)^\frac12.
\end{split}\]
But since $\mathcal V$, $D_\beta\nabla b$, and $D_\beta \nabla \phi$ all satisfies the Carleson measure condition, we conclude that $|T_{1221}| \lesssim r^d$ as desired.

We bounded each term derived from $T$ by either $r^d$ or $r^{d/2}T^{1/2}$, and consequently proved the claim \eqref{claimZ}. 
As was observed before, \eqref{GiCMz}, Proposition \ref{Main21}, and then Theorems \ref{Main2} and \ref{Main2a} follow. \ep

\end{document}